\documentclass[twoside]{article}
\usepackage[dvips]{epsfig}
\def\sqr#1#2{{\vcenter{\vbox{\hrule height.#2pt
              \hbox{\vrule width.#2pt height#1pt \kern#1pt \vrule width.#2pt}
              \hrule height.#2pt}}}}

\def\dbR{{\mathop{\rm l\negthinspace R}}}

\def\3n{\negthinspace \negthinspace \negthinspace }
\def\2n{\negthinspace \negthinspace }

\def\dbR{{\mathop{\rm l\negthinspace R}}}
\def\={\buildrel \triangle \over =}

\def\a{\alpha}

\def\d{\delta}
\def\e{\varepsilon}

\def\t{\times}

\def\o{\omega}

\def\ns{\noalign{\ss} }

\def\ds{\displaystyle}
\def\G{\Gamma}
\def\D{\Delta}

\def\Si{\Sigma}

\def\O{\Omega}

\def\cA{{\cal A}}

\def\cl{{\cal l}}
\def\no{\noindent}

\def\ss{\smallskip}
\def\ms{\medskip}

\def\q{\quad}
\def\qq{\qquad}
\def\hb{\hbox}

\def\pa{\partial}

\def\wt{\widetilde}
\def\cd{\cdot}

\def\dim{\hbox{\rm dim$\,$}}

\def\supp{\hbox{\rm supp$\,$}}

\def\cl{\overline}

\def\|{\Big |}
\def\({\Big (}
\def\){\Big )}
\def\[{\Big[}
\def\]{\Big]}
\def\be{\begin{equation}}
\def\bel{\begin{equation}\label}
\def\ee{\end{equation}}
\def\bt{\begin{theorem}}
\def\et{\end{theorem}}
\def\bc{\begin{corollary}}
\def\ec{\end{corollary}}
\def\bl{\begin{lemma}}
\def\el{\end{lemma}}
\def\bp{\begin{proposition}}
\def\ep{\end{proposition}}
\def\br{\begin{remark}}
\def\er{\end{remark}}
\def\ba{\begin{array}}
\def\ea{\end{array}}
\def\ed{\end{document}}

\newtheorem{lemma}{Lemma}[section]
\newtheorem{remark}{Remark}[section]

\newtheorem{theorem}{Theorem}[section]
\newtheorem{corollary}{Corollary}[section]

\newtheorem{proposition}{Proposition}[section]

\makeatletter
   
   \@addtoreset{equation}{section}
\makeatother

\begin{document}

\title
{\Large \bf \boldmath Remarks on the Controllability of Some
Quasilinear Equations\footnote{This work was supported by the NSF of
China under grants 10525105, 10831007 and 60821091, the Chunhui
program (State Education Ministry of China) under grant
Z007-1-61006, and the project MTM2008-03541 of the Spanish Ministry
of Science and Innovation. Part of this work was done when the
author visited Fudan University, with a financial support from
the ``French-Chinese Summer Institute on Applied Mathematics" (September 1-21, 2008).}}%

\author{\large Xu Zhang
    \\ \normalsize\emph{Academy of Mathematics and Systems
Sciences, }
    \\ \normalsize\emph{Chinese Academy of Sciences, Beijing 100190, China; and}
    \\ \normalsize\emph{Yangtze Center of Mathematics, }
    \\ \normalsize\emph{Sichuan University, Chengdu 610064,
China.}
    \\ \normalsize\emph{E-mail: xuzhang@amss.ac.cn}
}

\date{\vspace{-12mm}}

\maketitle


\begin{abstract}
In this Note, we review the main existing results, methods, and some
key open problems on the controllability of nonlinear hyperbolic and
parabolic equations. Especially, we describe our recent universal
approach to solve the local controllability problem of quasilinear
time-reversible evolution equations, which is based on a new
unbounded perturbation technique. It is also worthy to mention that
the technique we develop can also be applied to other problems for
quasilinear equations, say local existence,  stabilization, etc.
\end{abstract}

\section{Introduction}
Consider the following controlled evolution equation:
 \bel{1.0}
 \left\{
 \ba{ll}
 \displaystyle {d\over dt}y(t)=A(y(t))y(t)+Bu(t),\q t\in (0,T),\\
 \noalign{\ms}
 y(0)=y_0.
 \ea
 \right.
 \ee
Here, the time $T > 0$ is given, $y(t)\in Y$ is the state variable,
$u(t)\in U$ is the control variable, $y_0(\in Y)$ is the initial
state; $Y$ and $U$ are respectively the state space and control
space, both of which are some Hilbert space; $A(\cd) $ is a suitable
(nonlinear and usually unbounded) operator on $Y$, while the control
operator $B$ maps $U$ into $Y$. Many control problems for relevant
nonlinear Partial Differential Equations (PDEs, for short) enter
into this context. For instance, the quasilinear/semilinear
parabolic equation,  wave equation,  plate equation,  Schr\"odinger
equation, Maxwell equations, and  Lam\'e system, etc.

In this Note, we shall describe some existing methods, results and
main open problems on the controllability of these systems,
especially these for nonlinear hyperbolic and parabolic equations.

System (\ref{1.0}) is said to be exactly controllable in $Y$ at time
$T$ if for any $y_0,y_1\in Y$, there is a control $u\in L^2(0,T;U)$
such that the solution of system (\ref{1.0}) with this control
satisfies
  \bel{1.090}
 y(T)=y_1.
 \ee
When $\dim Y=\infty$ (We shall focus on this case later unless other
stated), sometimes one has to relax the requirement (\ref{1.090}),
and this leads to various notions and degrees of controllability:
approximate controllability, null controllability, etc. Note however
that, for time reversible system, the notion of exact
controllability is equivalent to that of null controllability.

Roughly speaking, the controllability problem for an evolution
equation consists in driving the state of the system (the solution
of the controlled equation under consideration) to a prescribed
final target state (exactly or in some approximate way) in finite
time. Problems of this type are common in science and engineering
and, in particular, they arise often in the context of flow control,
in the control of flexible structures appearing in flexible robots
and in large space structures, in quantum chemistry, etc.

The controllability theory  for finite dimensional linear systems
was introduced by R.E.~Kalman \cite{Kalman} at the very beginning of
the 1960s. Thereafter, many authors were devoted to develop it for
more general systems including infinite dimensional ones, and its
nonlinear and stochastic counterparts.

The controllability theory of PDEs depends strongly on its nature
and, in particular, on its time-reversibility properties. To some
extent, the study of controllability for linear PDEs is
well-developed although many challenging problems are still
unsolved. Classical references in this field are D.L.~Russell
\cite{Russell} and J.L.~Lions \cite{Lions}. Updated progress can be
found in a recent survey by E.~Zuazua (\cite{Zuazua}). Nevertheless,
much less are know for nonlinear controllability problems for PDEs
although several books on this topic are available, say J.M.~Coron
\cite{Coron}, A.V.~Fursikov $\&$ O.Yu.~Imanuvilov \cite{FI}, T.T.~Li
\cite{Li}, and X.~Zhang \cite{Zhang1}. Therefore, in this Note, we
concentrate on controllability problems for systems governed by
nonlinear PDEs.

The main result in this Note can be described as follows: Assume
that $(A(0),B)$ is exact controllable in $Y$. Then, under some
assumptions on the structure of $A(y)$ (for concrete problems, which
needs more regularity on the state space, say $\mathcal{D}(A(0)^k)$
for sufficiently large $k$), system (\ref{1.0}) is locally exact
controllable in $\mathcal{D}(A(0)^k)$.

The main approach that we employ to show the above controllability
result is a new perturbation technique. The point is that, the
perturbation is {\it unbounded} but {\it small}. Note however that
this approach does NOT work for the null controllability problem of
the time-irreversible systems, and therefore, one has to develop
different method to solve the local null controllability of
quasilinear parabolic equations.

For simplicity, in what follows, we consider mainly the case of
internal control, i.e. $B\in  \mathcal{L}(U,Y)$. Also, we will focus
on the local controllability of the quasilinear wave equation.
However, our approach is universal, and therefore, it can be
extended to other quasilinear PDEs, say quasilinear plate equation,
Schr\"odinger equation, Maxwell equations, and  Lam\'e system, etc.

On the other hand, we mention that the technique developed in this
Note can also be applied to other problems for quasilinear
equations. For example, stabilization problem for system (\ref{1.0})
(with small initial data) can be considered similarly. Indeed,
although there does not exist the same equivalence between exact
controllability and stabilization in the nonlinear setting, the
approaches to treat them can be employed each other.

The rest of this Note is organized as follows. In Section \ref{s2},
we review the robustness of the controllability in the setting of
Ordinal Differential Equations (ODEs, for short). In Section
\ref{s3}, we recall some known perturbation result on the exact
controllability of abstract evolution equations. Then, in Section
\ref{s4}, we show a new perturbation result on the exact
controllability of general evolution equations. Sections \ref{s5}
and \ref{s6} are addressed to present local controllability results
for multidimensional quasilinear hyperbolic equations and parabolic
equations, respectively. Finally, in Section \ref{s7}, we collect
some open problems, which seem to important in the field of
controllability of PDEs.

\section{Starting point: the case of ODEs}\label{s2}

Consider the following controlled ODE:
 \bel{00001.1}
 \left\{
 \ba{ll}
 \displaystyle {d\over dt}y=Ay+Bu,\qq t\in (0,T),\\
 \noalign{\ms}
 y(0)=y_0,
 \ea
 \right.
 \ee
where $A\in \dbR ^{n\t n}$ and  $B\in \dbR ^{n\t m}$. It is
well-known (\cite{Kalman}) that system (\ref{00001.1}) is exact
controllable in $(0,T)$ if and only if
  $$
  B^*e^{A^*t}x_0=0,\q \forall\; t\in (0,T)\Rightarrow x_0=0.
   $$
Note that this condition is also equivalen to the following Kalman
rank condition:
  \bel{llaa}
  rank (B,AB,A^2B,\cdots,A^{n-1}B)=n.
  \ee

From (\ref{llaa}), it is clear that if $(A, B)$ is exact
controllable, then there exists a small $\e=\e(A,B)>0$ such that
$(\wt A,\wt B)$ is still exact controllable provided that $||\wt
A-A||+||\wt B-B||<\e$. Therefore, the exact controllability of
system (\ref{00001.1}) is robust under small perturbation.

Because of the above robustness, the local exact controllability of
nonlinear OPEs is quite easy. Indeed, consider the following
controlled system:
 \bel{000--01.1}
 \left\{
 \ba{ll}
 \displaystyle {d\over dt}y=Ay+f(y)+Bu,\qq t\in (0,T),\\
 \noalign{\ms}
 y(0)=y_0,
 \ea
 \right.
 \ee
 with $f(\cd)\in C^1(\dbR ^n)$ and $f(y)=O(|y|^{1+\d})$ when $y$ is small, for some $\d>0$. The local exact
 controllability of system (\ref{000--01.1}) follows from a standard perturbation argument.

 However, the corresponding problem in PDE setting is much more complicated, as we shall see
 below.

\section{Known perturbation result on exact
controllability}\label{s3}

In this section, we recall some known perturbation result on the
exact controllability of abstract evolution equations. These results
are based on the following two tools:

\medskip

$\bullet$ {\bf Duality argument} (e.g. \cite{LY, Lions, Zhang1}): In
the linear setting (i.e., $A(y)\equiv A$ is independent of $y$ and
linear, and further $A$ generates an $C_0$-group
$\{e^{At}\}_{t\in\dbR }$ on $Y$), the null controllability of system
(\ref{1.0}) is equivalent to the following observability estimate:
 \bel{zx1}
 |e^{A^*T}z^*|_{Y^*}^2\le C\int_0^T|B^*e^{A^*s}z^*|_{U^*}^2ds,\qq\forall\, z^*\in
 Y^*,
 \ee
for some constant $C>0$.

\medskip

$\bullet$ {\bf Variation of constants formula}: In the setting of
semigroup, for a bounded perturbation $P\in {\cal L}(Y)$:
 \bel{zxe-1}
 \ba{ll}\ds
 e^{(A+P)t}x= e^{At}x+\int_0^te^{A(t-s)}Pxds,\qq\forall\, x\in
 Y.
 \ea
 \ee

\medskip

Combining (\ref{zx1}) and (\ref{zxe-1}), it is easy to establish the
following well-known (bounded) perturbation result on the exact
controllability:

\begin{theorem}\label{th3.1}
Assume that $A$ generates an $C_0$-group $\{e^{At}\}_{t\in\dbR }$ on
$Y$ and $B\in  \mathcal{L}(U,Y)$. If $(A,B)$ is exact controllable,
then so is $(A+P,B)$ provided that $||P||_{{\cal L}(Y)}$ is small
enough.
\end{theorem}

The above perturbation $P$ can also be time-dependent. In this case,
one needs the language of evolution system. In the sequel, for a
simple presentation, we consider only the time-independent case.

As a consequence of Theorem \ref{th3.1} and the standard fixed point
technique, one can easily deduce a local exact controllability
result for some semilinear equations, say the counterpart of system
(\ref{000--01.1}):
 \bel{0-01.1}
 \left\{
 \ba{ll}
 \displaystyle {d\over dt}z=Az+f(z)+Bv,\qq t\in (0,T),\\
 \noalign{\ms}
 z(0)=z_0.
 \ea
 \right.
 \ee
More  precisely, we have

\begin{corollary}\label{co3.1}
Assume that $A$ generates an $C_0$-group $\{e^{At}\}_{t\in\dbR }$ on
$Y$, $B\in  \mathcal{L}(U,Y)$, and $(A,B)$ is exact controllable. If
the nonlinearity $f(\cd):\; Y\to Y$ satisfies $f(\cd)\in C^1(Y)$
and, for some $\d>0$, $|f(z)|_Y=O(|z|_Y^{1+\d})$ as $|z|_Y\to 0$,
then system (\ref{0-01.1}) is locally exact controllable in $Y$.
\end{corollary}

Clearly, both the time reservability of the underlying system and
the variation of constants formula (\ref{zxe-1}) plays a key role in
the above perturbation-type results.

When the system is time-irreversible, the above perturbation
technique does not work. The typical example is the controlled heat
equation. In this case, one has to search for other robust method to
derive the desired controllability, say, Carleman estimate. We shall
consider this case in Section \ref{s6}.

When the perturbation operator $P$ is unbounded, formula
(\ref{zxe-1}) may fail to work, and in this case things become much
more delicate even for the semigroup theory itself. Nevertheless,
there do exist some special case, for which the perturbation $P$ is
unbounded but the above variation of constants formula still works
(in the usual sense), say when the semigroup $\{e^{At}\}_{t\ge 0}$
has some smooth effect. In this case, one can find some perturbation
result for exact controllability in S. Boulite, A. Idrissi and L.
Maniar \cite{BIM}, S. Hadd \cite{Hadd}, and H. Leiva \cite{Leiva}.
However, it does not seem that these perturbation results can be
adapted to solve the nonlinear controllability problems, especially
for quasilinear equations.

\section{A new perturbation result on exact
controllability}\label{s4}

In this section, we present a new perturbation result on the exact
controllability of general evolution equations. The idea is simple,
and the key point is that the generation of an $C_0$-semigroup
$\{e^{At}\}_{t\ge 0}$ is robust with respect to a small perturbation
of the same ``order" with respect to the generator $A$.

Stimulated by quasilinear problem, we consider the following small
perturbation of the same ``order":
 $$
 P=P_0A,
 $$
where $P_0\in {\cal L}(Y)$ and $||P_0||<1$. That is, the perturbed
operator reads: $(I+P_0)A$. It is easy to show that, if $A$
generated a contractive $C_0$-semigroup, then so is $(I+P_0)A$.
Indeed, it is obvious that  $(I+P_0)A$ is dissipative in $Y$ with
the new scalar product $((I+P_0)^{-1}\cd,\cd)$, which induces a
norm, equivalent to the original one. Nevertheless, we remark that
the variation of constants formula does not work for $e^{(I+P_0)At}$
for this general case.

Thanks to the above observation, a new perturbation result for exact
controllability is shown in \cite{Zhang2}, which reads as follows:

\begin{theorem}\label{th4.1}
Assume that $A$ generates an unitary group $\{e^{At}\}_{t\in\dbR }$
on $Y$ and $B\in  \mathcal{L}(U,Y)$. If $(A,B)$ is exact
controllable, then so is is $(A+P,B)\equiv ((I+P_0)A, B)$ provided
that $||P_0||_{{\cal L}(Y)}$ is small enough.
\end{theorem}

Since the variation of constants formula does not work for
$e^{(I+P_0)At}$, the above result can not be derived as Theorem
\ref{th3.1}. Instead, we need to use Laplace transform and some
elementary tools from complex analysis to prove the desired result.

The above simple yet useful perturbation-type controllability result
can be employed to treat the local controllability problems for
quasilinear evolution-type PDEs with time-reversibility, as we shall
see in the next section.

\section{Local exact controllability for multidimensional quasilinear hyperbolic
 equations}\label{s5}

This section is addressed to the local exact controllability of
quasilinear hyperbolic
 equations in any space dimensions.

To begin with, let us recall the related known controllability
results for controlled quasilinear hyperbolic equations. The problem
is well-understood in one space dimension. To the author's best
knowledge, the first paper in this direction is M.~Cirina
\cite{Cirina}. Recent rich results are available in T.T.~Li $\&$
B.P.~Rao \cite{LR}, T.T.~Li $\&$ B.Y.~Zhang \cite{LZ}, T.T.~Li $\&$
L.X.~Yu \cite{LY-1}, Z.Q.~Wang \cite{Wang}, and especially the above
mentioned book by T.T.~Li (\cite{Li}). As for the corresponding
controllability results in multi-space dimensions, we refer to
P.~F.~Yao (\cite{Yao}) and Y.~Zhou $\&$ Z.~Lei \cite{ZL}.

Let $\Omega$ be a bounded domain in $\dbR ^n$ with a sufficiently
smooth boundary $\Gamma$. Put $Q=(0,T)\times\Omega$ and
$\Sigma=(0,T)\times\Gamma$. Let $\omega$ be a nonempty open subset
of $\Omega$. We consider the following controlled quasilinear
hyperbolic equations:
   \bel{11oa5}
   \left\{\ba{ll}
   \ds z_{tt}-\sum_{i,j=1}^n\pa_{x_i}(a_{ij}(x) z_{x_j})\\
   \ns
   \ds\quad=G(t,x,z,\nabla_{t,x}z,\nabla_{t,x}^2z)+\phi_\o(x)u,\quad\quad&\hb{in }Q,\\
 \ns\ds z=0,\qq&\hb{in }\Si,\\
 \ns\ds  z(0)=z_0,z_t(0)=z_1,\qq&\hb{in }\O,
 \ea\right.\ee
where the coefficients $a_{ij}(\cdot)\in C^2(\cl{\O})$  ($i,
j=1,\cdots, n$) satisfy $a_{i j}=a_{j i}$, and for some constant
$\rho>0$,
$$
\displaystyle\sum^n_{i,j=1}a_{i j}(x)\xi_i\xi_j\geq \rho|\xi|^2, \ \
\forall\ (x, \xi)=(x,\xi_1,\cdots,\xi_n)\in \cl{\O}\times \dbR ^n;
$$
and following \cite{ZL}, the nonlinearity $G(\cd)$ is taken to be of
the form
 $$\ba{ll}\ds
G(t,x,\nabla_{t,x}z,\nabla_{t,x}^2z)\\
\ns
\ds=\sum_{i=1}^n\sum_{\a=0}^ng_{i\a}(t,x,\nabla_{t,x}z)\pa_{x_ix_\a}^2z+O(|u|^2+|\nabla_{t,x}z|^2),
 \ea
 $$
$g_{i\a}(t,x,0,0)=0$ and $x_0=t$; $\phi_\o$ is a nonnegative smooth
function defined on $\cl\O$ and satisfying $\ds\min_{x\in
\o}\phi_\o(x)>0$.

Denote by $\chi_\omega$ the characteristic function of $\omega$. We
need to introduce the following

\ms

\no {\bf Assumption (H)}: {\it Assume the linear hyperbolic equation
   \bel{oa5}
   \left\{\ba{ll}
   \ds y_{tt}-\sum_{i,j=1}^n\pa_{x_i}(a_{ij}(x) y_{x_j})=\chi_\o(x)u,\quad\quad&\hb{in }Q,\\
 \ns\ds y=0,\qq&\hb{in }\Si,\\
 \ns\ds  y(0)=y_0, y_t(0)=y_1,\qq&\hb{in }\O
 \ea\right.\ee
is exact controllable in $H_0^1(\O)\times L^2(\O)$}.

\ms

The following controllability result for quasilinear hyperbolic
equations is shown in \cite{Zhang2}:

\begin{theorem}\label{th5.1}
Let Assumption (H) hold. Then, for any $s>\frac{n}{2}+1$, system
(\ref{11oa5}) is {\it local exact controllable} in $(H^{s+1}(\O)\cap
H_0^1(\O))\times H^s(\O)$ (provided that some compatible conditions
are satisfied for the initial and final data).
\end{theorem}

Clearly, Theorem \ref{th5.1} covers the main results in \cite{Yao,
ZL}. The above result follows by combining our new perturbation
result for exact controllability , i.e. Theorem \ref{th4.1} and the
fixed point technique developed in \cite{ZL}.

\begin{remark}
The boundary control problem can be considered similarly although
the technique is a little more complicated.
\end{remark}

\begin{remark}
The key point of our approach is to reduce the local exact
controllability of quasilinear equations to the exact
controllability of the linear equation. This method is general and
simple. The disadvantage is that we can not construct the control
explicitly. Therefore, this approach does not replace the value of
\cite{ZL}, and the deep results for the corresponding $1-d$ problem,
obtained by T.~T.~Li and his collaborators, as mentioned before.
Especially, from the computational point of view, the later approach
might be more useful.
\end{remark}

We now return to Assumption (H), and review the known results and
unsolved problems for exact controllability of the linear hyperbolic
equation but we concentrate on the case of boundary control although
similar things can be said for the case of internal control.

Denote by $\cA$  the elliptic operator appeared in the first
equation of system (\ref{oa5}). We consider the following controlled
linear hyperbolic equation with a boundary controller:
 \bel{ollla5}
 \left\{\ba{ll}
   \ds y_{tt}+{\cal A} y=0,\quad\quad&\hb{in }Q,\\
 \ns\ds y=\chi_{\Si_0}u,\qq&\hb{in }\Si,\\
 \ns\ds  y(0)=y_0, y_t(0)=y_1,\qq&\hb{in }\O,
 \ea\right.\ee
where $\emptyset\not=\Si_0\subset\Si$ is the controller. It is easy
to show that, system (\ref{ollla5}) is exactly controllable in
$L^2(\O)\times H^{-1}(\O)$ at time $T$ by means of control $u\in
L^2(\Si_0)$ if and only if there is a constant $C>0$ such that
solutions of its dual system
 \bel{3-a5}\left\{\ba{ll}\ds w_{tt}+{\cal A} w=0,\qq&\hb{in }Q\\
 \ns\ds w=0,\qq&\hb{in }\Si\\
 \ns\ds  w(0)=w_0,w_t(0)=w_1,\qq&\hb{in }\O
 \ea\right.\ee
satisfies the following observability estimate:
 \bel{31-a5}
 \ba{ll}\ds
|w_0|_{H_0^1(\O)}^2+|w_1|_{L^2(\O)}^2\le
C\int_{\Si_0}\left|{\pa_{\cal A} w\over
\pa\nu}\right|^2d\Si_0,\\
 \ns
\ds\qq\qq\qq\qq\qq\forall\; (w_0,w_1)\in H_0^1(\O)\times L^2(\O).
\ea
 \ee

When $\cA=-\D$, $\Si_0=(0,T)\times\G_0$ with $\G_0$ to be a suitable
subset of $\pa\O$, L.F. Ho \cite{HO} establish (\ref{31-a5}) by
means of the classical Rellich-type multiplier. Later, K. Liu
\cite{Liu} gave a nice improvement for the case of internal control.
When $\cA$ is a general elliptic operator of second order, and
$\Si_0$ is a general (maybe non-cylinder) subset of $\Si$, J.L.
Lions \cite{Lions} posed an open problem on ``under which condition,
inequality (\ref{31-a5}) holds?". When $\Si_0=(0,T)\times\G_0$ is a
cylinder subset of $\Si$, Lions's problem is almost solved. In this
case, typical results are as follows:

\begin{enumerate}

\item[1)] Geometric
Optics Condition (GOC for short) introduced by C.~Bardos, G.~Lebeau
$\&$ J.~Rauch \cite{BLR}, which is a sufficient and (almost)
necessary condition for inequality (\ref{31-a5}) to hold. GOC is
perfect except the three disadvantage: One is that it needs
considerably high regularity on both the coefficients and $\pa\O$
(N. Burq \cite{Burq} gives some improvement in this respect); One is
that this condition is not easy to verify; The other is that the
observability constant derived from GOC is not explicit because it
involves the contradiction argument to absorb the undesired lower
order terms appeared in the observability estimate.

\item[2)] Rellich-type multiplier conditions introduced by L.F. Ho \cite{HO}, K. Liu
\cite{Liu}, A. Osses \cite{Osses}, etc., which require less smooth
conditions than GOC but they are not necessary conditions for
inequality (\ref{31-a5}) to hold.

\item[3)] There exist some other sufficient condition for
inequality (\ref{31-a5}) to hold, say the vector field condition by
A. Wyler \cite{Wyler}, and the curvature condition by P.F. Yao
(\cite{Yao0}. Later, it is shown by S.J.~Feng $\&$ D.X.~Feng
\cite{FF} that these two conditions are equivalent although they are
introduced through different tools.

\item[4)] Mixed tensor/vector field condition introduced by X. Zhang $\&$ E.
Zuazua \cite{ZZ}, which covers the conditions in 2) and 3).
\end{enumerate}

\begin{remark}
It is shown by L. Miller \cite{Miller} that when the data are
sufficiently smooth, the conditions in 2) and 3) are special cases
of GOC. Nevertheless, as far as I know, it is an unsolved problem on
the minimal assumption on data for GOC.
\end{remark}

When $\Si_0\not=(0,T)\times\G_0$, especially when it is NOT a
cylinder subset of $\Si$, there exist almost no nontrivial progress
on Lions's problem (which seems to be a challenging mathematical
problem), even for the simplest $1-d$ wave equation! The only
related results are as follows:

\begin{enumerate}

\item[a)] For $1-d$ wave
equation and $\Si_0=E\times\G_0$ with $E\subset (0,T)$ to be a
Lebesgue measurable set with positive measure, P. Martinez $\&$ J.
Vancostenoble \cite{MV} show that (\ref{31-a5}) holds.

\item[b)] G. Wang \cite{Wang0} obtains an interesting internal
observability estimate for the heat equation in multi-space
dimensions, where the observer is $E\times\o$ with $E$ being the
same as in the above case and $\o$ to be any nonempty open subset of
$\O$.
\end{enumerate}

\section{Local null controllability for quasilinear parabolic
 equations}\label{s6}

In this section, we consider the local exact controllability of
quasilinear parabolic equations in any space dimensions.

As mentioned before, the perturbation technique does not apply to
the time irreversible system, exactly the case of parabolic
equations. Therefore, one has to search for other robust method to
derive the desired null controllability, say, Carleman estimate even
if the perturbation to the null-controllable system is very small
(even in the linear setting!).

We consider the following controlled quasilinear parabolic system
\begin{eqnarray}\label{2}
\left\{
\begin{array}{lll}
&y_t-\displaystyle\sum^{n}_{i,j=1}(a_{i j}(y) y_{x_i})_{x_j}=\chi_\omega u  \ \ \ \ &\mbox{ in }Q,\\
&y=0 &\mbox{ on }\Sigma,\\
&y(0)=y_0 &\mbox{ in }\Omega,
\end{array}
\right.
\end{eqnarray}
where $a_{i j}(\cdot): \dbR \rightarrow \dbR $ are twice
continuously differentiable functions satisfying similar conditions
in the last section.

In the last decades, there are many papers devoted to the
controllability of linear and semilinear parabolic equations (see
e.g. \cite{FI, Zuazua} and the rich references therein). However, as
far as we know, nothing is known about the controllability of
quasilinear parabolic equations except for the case of one space
dimension. In \cite{Beceanu}, the author proves the local null
controllability of a $1-d$ quasilinear diffusion equation by means
of the Sobolev embedding relation
$L^\infty(0,T;H^1_0(\Omega))\subseteq L^\infty(Q)$, which is valid
only for one space dimension.

The following local null controllability result for a class of
considerably general multidimensional quasilinear parabolic
equations, system (\ref{2}), is shown in \cite{LZ}.

\medskip

\begin{theorem}\label{6t1}
There is a constant $\gamma>0$ such that, for any initial value
$y_0\in C^{2+\frac{1}{2}}(\overline{\Omega})$ satisfying
$|y_0|_{C^{2+\frac{1}{2}}(\overline{\Omega})}\leq \gamma$ and the
first order compatibility condition, one can find a control $u\in
C^{\frac{1}{2}, \frac{1}{4}}(\overline{Q})$ with $\supp
u\subseteq\omega\times[0,T]$ so that the solution $y$ of system
(\ref{2}) satisfies $y(T)=0$ in $\Omega$. Moreover,
 $$
|u|_{C^{\frac{1}{2}, \frac{1}{4}}(\overline{Q})}\leq C
e^{e^{CA}}|y_0|_{L^2(\Omega)},
 $$
where $A=\sum\limits^n_{i,j=1}\left(1+\sup\limits_{|s|\leq 1}|a_{i
j}(s)|^2+\sup\limits_{|s|\leq 1}|a_{i j}'(s)|^2\right)$, and $C$
depends only on $\rho$, $n$, $\Omega$ and $T$.
\end{theorem}

The key point in the proof of Theorem \ref{6t1} is to improve the
regularity of the control function for smooth data, which is a
consequence of a new observability inequality for linear parabolic
equations with an explicit estimate on the observability constant in
terms of the $C^1$-norm of the coefficients in the principle
operator. The later is based on a new global Carleman estimate for
the parabolic operator.

\section{Open problems}\label{s7}

Although great progress have been made on the controllability theory
of PDEs, the field is still full of open problems. In some sense,
the linear theory is well-understood and there exist extensive works
on the controllability of linear PDEs. But, still, even for the
linear setting, some fundamental problems remain to be solved, as we
shall explain later. The controllability theory of nonlinear system
originated in the middle of 1960s but the progress is very slow.
Similar to other nonlinear problems, controllability of infinite
dimensional nonlinear system is usually very difficult. Due to the
underlying properties of the equation, the progress of the exact
controllability theory for nonlinear hyperbolic equations is even
slower. Nevertheless, nonlinear problems are not always difficult
than linear ones. Indeed, as we have shown in Theorem \ref{th5.1},
local exact controllability of quasilinear hyperbolic equations is a
consequence of the exact controllability of linear hyperbolic
equations. One may then ask such a question: ``How to judge a
nonlinear result is good or not?" To the author's opinion, except
for some famous unsolved problem, the point is either ``whether the
result is optimal or not in some nontrivial sense?", or ``whether
some new phenomenon is discovered or not?".

From the above ``criteria", our result on the local exact
controllability of quasilinear hyperbolic equations is not good at
all. Indeed, there is no evidence to show that the result is
optimal. Therefore,
\begin{enumerate}

\item[] {\it How to establish the ``optimal" local exact controllability \\ result
for quasilinear equations?}

\end{enumerate}
\no is one of the most challenging problems in the field of control
of PDEs. As we shall see below, this problem is also highly
nontrivial even in the semilinear setting!

We now review the exact controllability for the following semilinear
hyperbolic equations:
  \bel{1-1-a5}\left\{\ba{ll}\ds z_{tt}+\cA z=f(z)+\chi_\o(x)u(t,x),\qq&\hb{in }Q,\\
 \ns\ds z=0,\qq&\hb{in }\Si,\\
 \ns\ds  z(0)=z_0,z_t(0)=z_1,\qq&\hb{in }\O.
 \ea\right.
 \ee
For some very general nonlinearity $f(\cd)$ and a suitable
controller $\o$, E. Zuazua \cite{Zu0} obtains the local exact
controllability for system (\ref{1-1-a5}). Recently, B. Dehman \& G.
Lebeau \cite{DehmanL} gave a significant improvement. However, as
far as I know, no optimality on the controllability results are
analyzed in these works, which seems also to be a challenging
problem.

\begin{remark}
The possible optimality on the local exact controllability
for semilinear equations should be strongly related to PDEs with
lower regularity dada. This is a very rapid developing field in
recent years.
\end{remark}

\begin{remark}
There exists big difference between the controllability problems and
pure PDEs problems. Indeed, the exact controllability problem for
the system
  \bel{1-1-ba5}\left\{\ba{ll}\ds z_{tt}+\cA z=f(z_t)+\chi_\o(x)u(t,x),\qq&\hb{in }Q\\
 \ns\ds z=0,\qq&\hb{in }\Si\\
 \ns\ds z(0)=z_0,z_t(0)=z_1,\qq&\hb{in }\O
 \ea\right.\ee
in the natural energy space $H^1_0(\O)\times L^2(\O)$ is not clear
even if $f(\cd)$ is global Lipchtiz continuous. But, of course, the
well-posedness of the corresponding pure PDE problem (i.e. the
control $u\equiv 0$) is trivial.
\end{remark}

Global exact controllability for semilinear equations is generally a
very difficult problem. We refer to \cite{Zhang1} for known global
controllability results for the semilinear hyperbolic equation when
the nonlinearity is global Lipschitz continuous. For system
(\ref{1-1-a5}), if the nonlinearity $f(\cd)$ grows too fast, say
 \begin{equation}\label{1w.5}
 \lim_{|s|\to\infty}
 |f(s)||s|^{-1}\log^{-r}|s|=0,\qq r>2,
 \end{equation}
the solution may blowup, and therefore, global exact controllability
is impossible in this case. Recently, based on X.~Fu, J.~Yong $\&$
X.~Zhang \cite{FYZ} and V.Z.~Meshkov \cite{Me}, T.~Duyckaerts,
X.~Zhang $\&$ E.~Zuazua \cite{DZZ} showed that, if
 \begin{equation}\label{1eeew.5}
 \lim_{|s|\to\infty}
 |f(s)||s|^{-1}\log^{-r}|s|=0,\qq r<3/2,
 \end{equation}
then system (\ref{1-1-a5}) is globally exact controllable. Moreover,
it is also shown that the above index ``$3/2$" is optimal in some
sense (i.e., wether the linearization argument works or not) when
$n\ge 2$. (But this number is not optimal in $1-d$).

\begin{remark}
The same ``$3/2$"-phenomenon happens also for parabolic equations
when $n\ge 2$. Surprisingly, the $1-d$ problem is unsolved. That is,
it is not clear whether the index ``$3/2$" is optimal or not in
$1-d$! This means, sometimes, the $1-d$ problem is difficult than
the multidimensional ones.

\end{remark}

\begin{remark}
 Note that for the pure PDE problems, the same phenomenon
described above does not happen. This indicates that the study of
the controllability problem for nonlinear PDEs has some independent
interest, which is far from a sub-PDE-problem.
\end{remark}

\begin{remark}
Another strongly related longstanding unsolved problem is the exact
controllability of the linear time- and space-dependent hyperbolic
equation under the GOC. It seems that, this needs to combine
cleverly the tool from micro-local analysis and the technique of
Carleman estimate. But nobody knows how to do it.
\end{remark}

To end this Note, we list the following further open problems.

\begin{enumerate}

\item[$\bullet$] {\bf Controllability of the coupled and/or higher order systems
by using minimal number of controls}. As shown in X. Zhang $\&$ E.
Zuazua \cite{ZZ1}, the study of the related controllability problem
is surprisingly complicated and highly nontrivial even for the
systems in one space dimension!

\item[$\bullet$]{\bf Constrained controllability}. As shown in K.D.~Phung,
G.~Wang $\&$ X. Zhang \cite{PWZ}, the problem is unexpected
difficult even for the simplest $1-d$ wave equation and heat
equation.

\item[$\bullet$]{\bf Controllability of parabolic PDEs with memory, or retard argument
and/or other nonlocal terms}. Consider the following controlled heat
equations with a memory term:
  $$
  \left\{\ba{ll}\ds z_{t}-\D z=\int_0^ta(s,x)z(s)ds+\chi_\o(x)u,\qq&\hb{in }Q,\\
 \ns\ds z=0,\qq&\hb{in }\Si,\\
 \ns\ds  z(0)=z_0,\qq&\hb{in }\O.
 \ea\right.
  $$
The PDE problem itself is not difficult. But, as far as I know, the
controllability problem for the above equation is unsolved even if
the memory kernel $a(\cd,\cd)$ is small!

\item[$\bullet$] {\bf  Controllability/observability of stochastic PDEs}. There exists only
very few nontrivial results, say \cite{TZ, Zhang3} and the reference
cited therein. I believe this is a very hopeful direction for the
control of PDEs in the near future.

\item[$\bullet$] {\bf  Controllability of PDEs in non-reflexive space}. There
exists almost no nontrivial results in this direction!

\item[$\bullet$] {\bf Other types of controllability}. Different notions of controllability, say, periodic
controllability, may lead to new and interesting problems for PDEs.

\end{enumerate}


\end{document}